\documentclass[11pt]{article}
\usepackage{amsmath,amssymb,amsbsy,amsfonts,amsthm,latexsym,
            amsopn,amstext,amsxtra,euscript,amscd,amsthm}

\def\\{\cr}
\def\({\left(}
\def\){\right)}
\def\[{\left[}
\def\]{\right]}
\def\<{\langle}
\def\>{\rangle}

\begin{document}

\title{Corrigendum to ``Cullen numbers with the Lehmer Property"}

\author{{\sc Jos\'e~Mar\'{\i}a~Grau~Ribas}\\
Departamento de M\'atematicas\\
Universidad de Oviedo\\
 Avda. Calvo Sotelo, s/n, 33007 Oviedo, Spain\\
 {\tt grau@uniovi.es}
 \and
{\sc Florian~Luca}\\ 
Instituto de Matem{\'a}ticas,\\
Universidad Nacional Autonoma de M{\'e}xico,\\
C.P. 58089, Morelia, Michoac{\'a}n, M{\'e}xico\\
{\tt fluca@matmor.unam.mx}
}

\date{\today}

\pagenumbering{arabic}

\maketitle

\begin{abstract}
In this note, we correct an oversight from the paper \cite{GL} mentioned in the title.
\end{abstract}

There is an error on Page 131 of the paper \cite{GL} mentioned in the title in justifying that the expression $A$ is nonzero. After the sentence
``Also, since $m_p$ divides $n_1$, it follows that $u\le w$" on Page 131 in \cite{GL}, the argument continues in the following way. The case when $\rho=1$ implies $n_1=1$ and leads to the conclusion
that all prime factors of $C_n$ are Fermat primes, and this instance has been dealt with on Page 131 in \cite{GL}. Thus, we may assume that $\rho\ge 3$. The relation
$$
(2^{\alpha}\rho^w+\alpha)u=wn_p
$$
shows that $u\mid n_p$. Thus, 
$$
p=m_p 2^{n_p}+1=\rho^u 2^{n_p}+1=X^{u}+1,
$$
where $X=\rho 2^{n_p/u}$ is an integer. If $u>1$, the above expression has $X+1$ as a proper divisor $>1$ (because $u$ is odd), which is impossible since $p$ is prime. Thus, $u=1$. If $w=1$, we first get  that $m_p=n_1=\rho$, and then that $n_p=\alpha+2^{\alpha}\rho=\alpha+n$, so $p=C_n$, which is not allowed. Otherwise, $w\ge 3$,  $n_1=\rho^{w}$ and $p=\rho 2^{(\alpha+n)/w}+1=(n2^n)^{1/w}+1$. We now show that there is at most one prime $p$ with the above property. Indeed, assume that there are two of them $p_1$ and $p_2$, corresponding to $w_1<w_2$. Thus, 
$n_1=\rho_1^{w_1}=\rho_2^{w_2}$ and both $w_1$ and $w_2$ divide $n+\alpha$. Let $W={\text{\rm lcm}}[w_1,w_2]$. Then 
$n_1=\rho_0^W$ for some positive integer $\rho_0$. Furthermore, writing $W=w_1\lambda$, we have that $\lambda>1$, and $\rho_0^\lambda=\rho_1$. Hence,
$$
p_1=\rho_1 2^{(\alpha+n)/w_1}+1=Y^{\lambda}+1,
$$
where $Y=\rho_0 2^{(\alpha+n)/W}$ is an integer. This is false since $\lambda>1$ is odd, therefore the above expression $Y^{\lambda}+1$ has  $Y+1$ as a proper 
divisor $>1$, contradicting the fact that $p_1$ is prime. Hence, if $A$ is zero for some $p$, then $p$ is unique. Further, in this case $n_1=\rho^w$ and $p=(n2^n)^{1/w}+1\le 
(n2^n)^{1/3}+1$. 

The remaining of the argument from the paper  \cite{GL} shows that the  expression $A$ is nonzero for all other primes $q$ of $C_n$, so all prime factors 
$q$ of $C_n$ satisfy inequality (5) in the paper \cite{GL} with at most one exception, say $p$, which satisfies the inequality $p\le (n2^n)^{1/3}
+1$. Hence, instead of the inequality from Line 2 of Page 132 in \cite{GL}, we get that
$$
C_n<((n2^n)^{1/3}+1) 2^{6(k-1)(n\log n)^{1/2}},
$$
giving
$$
2^{6(k-1)(n\log n)^{1/2}}>\frac{n2^n}{(n2^n)^{1/3}+1}>2^{2n/3},
$$
where the right--most inequality above holds for all $n\ge 3$. This leads to a slightly worse inequality than the inequality (6) in the paper \cite{GL}, namely
\begin{equation}
\label{eq:lowk}
k>1+\frac{n^{1/2}}{9(\log n)^{1/2}}.
\end{equation}
Note that inequality (6) from the paper \cite{GL} still holds whenever $A\ne 0$ for all primes $p$ dividing $n$, and in particular  for all $n$ except maybe when $n_1=\rho^w$ for some $\rho\ge 3$ and $w\ge 3$. So, from now on, we shall treat only the case when $n_1=\rho^w$. 
Comparing estimate (3) in the paper \cite{GL} with \eqref{eq:lowk} leads to 
\begin{equation}
\label{eq:*}
\frac{n^{1/2}}{9(\log n)^{1/2}}<2.4\log n,
\end{equation}
which implies that $n<1.4\times 10^6$.  We now lower the bound in a way similar to the calculation on Page 132 in \cite{GL}. Namely, first if $2^{2^{\gamma}}+1$ is a Fermat prime factor of $C_n$, then $\gamma\le 20$, so $\gamma\in \{0,1,2,3,4\}$. Furthermore, $\log n/\log 3\le 12.9$, therefore $k\le 5+12=17$. Now inequality \eqref{eq:lowk} shows that
$$
\frac{n^{1/2}}{9(\log n)^{1/2}}<16,
$$
giving $n<260,000$. But then $\log n/\log 3\le 11.4$, giving $k\le 16$. Also, if $n$ is not a multiple of $3$, then the number of prime factors $p$ of $C_n$ with 
$m_p>1$ is at most $\log 260,000/\log 5<7.8.$ Thus, $C_n$ can have at most $5+7=12$ distinct prime factors, contradicting the result of Cohen and Hagis \cite{CoHa}. Hence, $3\mid n$ showing that $3$ does not divide $C_n$.  Thus, $k\le 15$, so
$$
\frac{n^{1/2}}{9(\log n)^{1/2}}<14,
$$
giving $n<200,000$. Also, $n$ cannot be divisible by a prime $q\ge 5$, for otherwise, since $n_1=\rho^w$ for some $w\ge 3$, we would get that the number of prime factors $p$ of $C_n$ with $m_p>1$ is at most $3+\log(200,000/q^3)/\log 3<9.8$, so $k\le 9+4=13$, contradicting again the result of Cohen and Harris. Hence, $n=2^{\alpha}\cdot 3^{\beta}$ 
and the proof finishes as in the paper \cite{GL} after formula (7).

\medskip

{\bf Acknowledgements.} We thank Dae Jun Kim for pointing out the oversight.

\end{document}